\documentclass{amsart}

\usepackage{fancyhdr}
\usepackage{lastpage}
\usepackage{stmaryrd,yhmath}

\newcommand{\bdis}{\begin{displaymath}}
\newcommand{\edis}{\end{displaymath}}
\newcommand{\be}{\begin{equation}}
\newcommand{\ee}{\end{equation}}

\newcommand{\mcal}{\mathcal}

\newcommand{\vt}{\vartheta}

\newcommand{\zf}{\zeta\left(\frac{1}{2}+it\right)}

\newtheorem{lemma}[]{Lemma}

\theoremstyle{definition}

\newtheorem{cor}[]{Corollary}

\theoremstyle{remark}
\newtheorem{remark}[]{Remark}

\newtheorem*{mydef1}{{\bf Theorem}}

\numberwithin{equation}{section}

\begin{document}

\title{A minimal integral of the Riemann $\Xi$-function}

\author{Jan Moser}

\address{Department of Mathematical Analysis and Numerical Mathematics, Comenius University, Mlynska Dolina M105, 842 48 Bratislava, SLOVAKIA}

\email{jan.mozer@fmph.uniba.sk}


\begin{abstract}
In this paper we obtain an equilibrium sequence $\{ \omega_n\}$ for which the following holds true: the areas (measures) of the figures corresponding
to the positive and negative parts, respectively, of the graph of the function $\Xi(t),\ t\in[\omega_n,\omega_{n+1}]$ are equal. \\

\noindent Dedicated to the 500th anniversary of rabbi L\" ow.
\end{abstract}

\maketitle

\section{The result}

\subsection{}

Hardy proved in 1914 the following fundamental theorem: the function $\zf$ has infinitely many real zeros (see \cite{1}). To prove this Hardy used
the following complicated formula

\be \label{1.1}
\int_0^\infty\frac{\Xi(t)}{t^2+\frac 14}t^{2n}\cosh(\mu t){\rm d}t\sim (-1)^n\frac{\pi}{2^{2n}}\cos\frac \pi 8 ,\ \mu\to\frac \pi 4-0
\ee
where
\bdis
\Xi(z)=\xi\left(\frac 12+iz\right),\ \xi(s)=\frac{s(s-1)}{2}\pi^{-s/2}\Gamma\left(\frac{s}{2}\right)\zeta(s),\ s=\sigma+it .
\edis

Let us note that (\ref{1.1}) follows from one well-known Ramanujan's formula (see \cite{5}, p. 36). \\

In this direction we obtain the following

\begin{mydef1}
There is an increasing sequence $\{ \omega_n\}_{n=1}^\infty$ such that
\be \label{1.2}
\int_{\omega_n}^\infty \Xi(t){\rm d}t=0 ,\ n=1,2,\dots  ,
\ee
and
\be \label{1.3}
\omega_{n+1}-\omega_n<\omega_n^{\frac 16+\epsilon} .
\ee
\end{mydef1}

\begin{remark}
The integral (\ref{1.2}) is to be named \emph{the minimal integral} of the Riemann $\Xi(t)$-function (comp. (\ref{1.1})).
\end{remark}

\subsection{}

From (\ref{1.2}) follows immediately

\begin{cor}
\be \label{1.4}
\int_{[\omega_n,\infty)^+}\Xi(t){\rm d}t=-\int_{[\omega_n,\infty)^-}\Xi(t){\rm d}t,\quad n=1,2,\dots
\ee
where
\bdis
[\omega_n,\infty)^+=\{\ t:\ \Xi(t)>0,\ t\in [\omega_n,\infty)\},\ \dots  .
\edis
\end{cor}

\begin{remark}
The global law of the exact equality of the areas (measures) of the figures corresponding (by the usual way) to the positive and negative parts, respectively,
of the graph of the function $\Xi(t),\ t\in[\omega_n,\infty)$ is expressed by the formula (\ref{1.4}).
\end{remark}

Next, from (\ref{1.2}) we obtain

\begin{cor}
\be\label{1.5}
\int_{\omega_n}^{\omega_{n+1}}\Xi(t){\rm d}t=0,\quad n=1,2,\dots \ ,
\ee
and consequently
\be \label{1.6}
\int_{[\omega_n,\omega_{n+1}]^+}\Xi(t){\rm d}t=-\int_{[\omega_n,\omega_{n+1}]^-}\Xi(t){\rm d}t .
\ee
\end{cor}

\begin{remark}
The local law of the exact equality of the areas (measures) of the figures corresponding (by the usual way) to the positive and negative parts, respectively,
of the graph of the function $\Xi(t),\ t\in[\omega_n,\omega_{n+1}]$ is expressed by the formula (\ref{1.6}).
\end{remark}

We obtain, by making use of the mean-value theorem in the formula (\ref{1.5}),
\bdis
\Xi(c)=0,\ c\in (\omega_n,\omega_{n+1}) \ \Rightarrow \ c=\gamma:\ \zeta\left(\frac 12+\i\gamma\right)=0 ,
\edis
i.e. we have
\begin{cor}
The following canonical property
\be \label{1.7}
\{ c\}_n=\{ \gamma\}_n,\ \{ c\}_n\subset (\omega_n,\omega_{n+1}) ,
\ee
holds true, i.e. the set $\{ c\}_n$ of the mean-value-points of the function $\Xi(t),\ t\in [\omega_n,\omega_{n+1}]$ is identical with the
set of zeros $\{ \gamma\}_n$ of this function ($\Xi(\gamma)=0\ \Leftrightarrow\ \zeta(\frac 12+i\gamma)=0$).
\end{cor}

\begin{remark}
By means of (\ref{1.4}), (\ref{1.6}), (\ref{1.7}) we have named the sequence $\{\omega_n\}_n$ as \emph{the equilibrium sequence}.
\end{remark}

\begin{remark}
Let us remind explicitly that the proof of the minimal integral (\ref{1.2}) is, at the same time, the new kind of the proof of the mentioned
Hardy's theorem (comp. (\ref{1.7})).
\end{remark}

\begin{remark}
The small improvements of the exponent $\frac 16$ in (\ref{1.3}) are irrelevant. In this direction, see our discussion connected with the
I.M. Vinogradov's scepticism on possibilities of the method of trigonometric sums (see \cite{4}).
\end{remark}

\section{The main formula}

\subsection{}

Let us remind the Riemann-Siegel formula
\be \label{2.1}
Z(t)=2\sum_{n\leq x(t)}\frac{1}{\sqrt{n}}\cos\{\vartheta(t)-t\ln n\}+\mcal{O}(t^{-1/4}),\ x(t)=\sqrt{\frac{t}{2\pi}} ,
\ee
(see \cite{5}, pp. 79. 221) where
\be \label{2.2}
\begin{split}
& Z(t)=e^{i\vartheta(t)}\zf,\ \vartheta(t)=-\frac 12t\ln\pi+\text{Im}\left\{\ln\Gamma\left(\frac 14+i\frac t2\right)\right\}= \\
&=\frac 12 t\ln\frac{t}{2\pi}-\frac 12 t-\frac 18\pi+\mcal{O}(t^{-1}),\ \vartheta'(t)=\frac 12\ln\frac{t}{2\pi}+\mcal{O}(t^{-1}),
\ \vartheta''(t)\sim\frac{1}{2t} ,
\end{split}
\ee
and the Gram's sequence $\{ t_\nu\}$ is defined by the equation (see \cite{5})
\be \label{2.3}
\vartheta(t_\nu)=\pi\nu,\quad \nu\geq \nu_0\geq 1 .
\ee

\subsection{}

Let

\be\label{2.4}
\begin{split}
& \Phi_1(T)=\int_T^\infty e^{-\alpha t}t^\beta Z(t){\rm d}t=\lim_{T'\to\infty}\int_T^{T'} e^{-\alpha t}t^\beta Z(t){\rm d}t=\lim_{T'\to\infty}S(T,T') , \\
& \alpha=\frac{\pi}{4},\ \beta=\frac{7}{4},\ Z(t)=\mcal{O}(t^{1/4}),\ T\geq T_0>0 ,
\end{split}
\ee
where $T_0$ is a sufficiently big. By the formula (\ref{2.1}) we have
\be \label{2.5}
S(T,T')=2S_1+2S_2+2S_3+Q
\ee
where
\be \label{2.6}
\begin{split}
& S_1=\int_T^{T'}e^{-\alpha t}t^\beta\cos\{\vartheta(t)\}{\rm d}t,\quad Q=\mcal{O}\left(\int_T^\infty e^{-\alpha t}t^{\beta-\frac 14}{\rm d}t\right), \\
& S_2=\sum_{2\leq n\leq\sqrt{\frac{T}{2\pi}} }\frac{1}{\sqrt{n}}\int_T^{T'}e^{-\alpha t}t^\beta\cos\{\vartheta(t)-t\ln n\}{\rm d}t , \\
& S_3=\sum_{\sqrt{\frac{T}{2\pi}}<n\leq\sqrt{\frac{T'}{2\pi}} }\frac{1}{\sqrt{n}}\int_{T_1=2\pi n^2}^{T'}e^{-\alpha t}t^\beta\cos\{\vartheta(t)-t\ln n\}{\rm d}t
\end{split}
\ee
by using the formula
\bdis
\int_T^{T'}\sum_{n\leq \sqrt{\frac{t}{2\pi}}}=\sum_{n\leq \sqrt{\frac{T'}{2\pi}}}\int_{T_1}^{T'},\quad T_1=\max\{ T,2\pi n^2\} .
\edis

\subsection{}

We obtain the following formula for the integral in $S_2$
\be \label{2.7}
\begin{split}
& \int_T^{T'}e^{-\alpha t}t^\beta\cos\{\vartheta(t)-t\ln n\}{\rm d}t=\frac{1}{\alpha}e^{-\alpha T}T^\beta\cos\{\vartheta(T)-T\ln n\}- \\
& -\frac{1}{\alpha^2}e^{-\alpha T}T^\beta\{\vartheta'(T)-\ln n\}\sin\{\vartheta(T)-T\ln n\} - \\
& -\frac{1}{\alpha^2}\int_T^{T'}e^{-\alpha t}t^\beta\{\vartheta'(t)-\ln n\}^2\cos\{\vartheta(t)-t\ln n\}{\rm d}t- \\
& -\frac{2\beta}{\alpha^2}\int_T^{T'}e^{-\alpha t}t^{\beta-1}\{\vartheta'(t)-\ln n\}\sin\{\vartheta(t)-t\ln n\}{\rm d}t- \\
& -\frac{1}{\alpha^2}\int_T^{T'}e^{-\alpha t}t^\beta\vartheta''(t)\sin\{\vartheta(t)-t\ln n\}{\rm d}t+ \\
& +\frac{\beta(\beta-1)}{\alpha^2}\int_T^{T'}e^{-\alpha t}t^{\beta-1}\cos\{\vartheta(t)-t\ln n\}{\rm d}t+ \\
& + \mcal{O}(e^{-\alpha T}T^{\beta-1})+\mcal{O}(e^{-\alpha T'}(T')^\beta) .
\end{split}
\ee

\subsection{}

Since (see (\ref{2.2}))

\be \label{2.8}
\begin{split}
& \{\vartheta'(t)-\ln n\}^2=\{\vartheta'(T)-\ln n\}^2+\mcal{O}\left\{ (t-T)\frac{\ln T}{T}\right\};\ \frac{\ln x}{x}<\frac{\ln T}{T},\ x>T ,
\end{split}
\ee
and
\be \label{2.9}
\begin{split}
& \int_T^{T'}e^{-\alpha t}t^\beta (t-T){\rm d}t= \\
& = -\frac 1\alpha e^{-\alpha T'}T'^\beta(T'-T)+\frac{\beta+1}{\alpha}\int_T^{T'}e^{-\alpha t}t^\beta{\rm d}t-\frac \beta\alpha
\int_T^{T'}e^{-\alpha t}t^{\beta-1}{\rm d}t= \\
& \dots =\mcal{O}\{ e^{-\alpha T'}T'^{\beta+1}\}+\mcal{O}(e^{-\alpha T}T^\beta)+\mcal{O}\left\{\int_T^{T'}e^{-\alpha t}t^{\beta-2}{\rm d}t\right\}= \\
& =\mcal{O}\{ e^{-\alpha T'}T'^{\beta+1}\}+\mcal{O}(e^{-\alpha T}T^\beta);\ \beta-2<0 ,
\end{split}
\ee
then ($\frac{\ln T}{T}<1,\ T>e$)
\be \label{2.10}
\frac{\ln T}{T}\int_T^{T'}e^{-\alpha t}t^\beta (t-T){\rm d}t=\mcal{O}\left( e^{-\alpha T}T^\beta\frac{\ln T}{T}\right)+
\mcal{O}\left\{ e^{-\alpha T'}T'^{\beta+1}\right\} .
\ee
Hence, for the first integral on the right-hand side of (\ref{2.7}) we obtain (see (\ref{2.8}), (\ref{2.10}))

\be \label{2.11}
\begin{split}
& \int_T^{T'}e^{-\alpha t}t^\beta\{\vartheta'(t)-\ln n\}^2\cos\{\vartheta(t)-t\ln n\}{\rm d}t= \\
& =\{\vt'(T)-\ln n\}^2\int_T^{T'}e^{-\alpha t}t^\beta\cos\{\vartheta(t)-t\ln n\}{\rm d}t+\mcal{O}\left( e^{-\alpha T}T^\beta\frac{\ln T}{T}\right)+ \\
& + \mcal{O}\{ e^{-\alpha T'}T'^{\beta+1}\} .
\end{split}
\ee

\subsection{}

For all the remaining integrals we obtain by the same method the estimate
\be \label{2.12}
\mcal{O}(e^{-\alpha T}T^{\beta-1}\ln T) .
\ee
Thus, from (\ref{2.7}) by (\ref{2.11}), (\ref{2.12}) we obtain the following formula (see $S_2$ in (\ref{2.6}))
\be \label{2.13}
\begin{split}
& \int_T^{T'} e^{-\alpha t}t^\beta\cos\{\vt(t)-t\ln n\}{\rm d}t= \\
& = \frac 1\alpha\frac{e^{-\alpha T}T^\beta}{1+\frac{1}{\alpha^2}\{\vt'(T)-\ln n\}^2}
\left[\cos\{\vt(T)-T\ln n\}-\frac{}{} \right. \\
& \left.-\frac{1}{\alpha}\{\vt'(T)-\ln n\}\sin\{\vt(T)-T\ln n\} \right]+\mcal{O}\left( e^{-\alpha T}T^\beta\frac{\ln T}{T}\right)+ \\
& +\mcal{O}\{ e^{-\alpha T'}T'^{\beta+1}\} ,
\end{split}
\ee
and consequently, for the integral $S_1$ (see (\ref{2.6})) we have
\be\label{2.14}
\begin{split}
& \int_T^{T'}e^{-\alpha t}t^\beta\cos\{\vt(t)\}{\rm d}t=\frac 1\alpha\frac{e^{-\alpha T}T^\beta}{1+\frac{1}{\alpha^2}\{\vt'(T)\}^2}
\left[\cos\{\vt(T)\}-\frac{}{} \right. \\
& \left. -\frac 1\alpha\vt'(T)\sin\{\vt(T)\}\right]+\mcal{O}\left( e^{-\alpha T}T^\beta\frac{\ln T}{T}\right)
+\mcal{O}\{ e^{-\alpha T'}T'^{\beta+1}\} .
\end{split}
\ee

\subsection{}

Since for the integral in $S_3$ (see (\ref{2.6}), comp. (\ref{2.9})) the estimate
\bdis
\int_{T_1}^{T'}e^{-\alpha t}t^\beta\cos\{\vt(t)-t\ln n\}{\rm d}t=\mcal{O}(e^{-\alpha T_1}T_1^\beta)=\mcal{O}(e^{-2\pi \alpha n^2}n^{2\beta})
\edis
holds true, then
\be \label{2.15}
e^{\alpha T}T^{-\beta}S_3=\mcal{O}
\left\{\sum_{\sqrt{\frac{T}{2\pi}}<n\leq\sqrt{\frac{T'}{2\pi}}}\frac{1}{\sqrt{n}}\left(\frac{n^2}{T}\right)^\beta
e^{\alpha(-2\pi n^2+T)}\right\} .
\ee
Putting in (\ref{2.15}) $n=n_0+k$ where
\bdis
n_0-1\leq\sqrt{\frac{T}{2\pi}}<n_0;\quad -2\pi n_0^2+T<0
\edis
we continue to manipulate with this estimate by the following way
\be \label{2.16}
\begin{split}
&=\mcal{O}
\left\{\sum_{\sqrt{\frac{T}{2\pi}}<n_0+k\leq\sqrt{\frac{T'}{2\pi}}}\frac{1}{\sqrt{n_0+k}}\left(\frac{n_0+k}{n_0-1}\right)^{2\beta}
e^{\alpha(-2\pi (n_0+k)^2+T)}\right\}= \\
& =\mcal{O}
\left\{\frac{1}{\sqrt{n_0}}\left[\left(\frac{n_0}{n_0-1}\right)^{2\beta}+\sum_{k=1}^\infty\left(\frac{n_0+k}{n_0-1}\right)^{2\beta}
e^{-2\pi \alpha(k^2+2n_0k)}\right]\right\}=\\
& = \mcal{O}\left(\frac{1}{\sqrt{n_0}}\right)=\mcal{O}(T^{-1/4}) \ \Rightarrow \
S_3=\mcal{O}(e^{-\alpha T}T^{\beta-\frac 14}) .
\end{split}
\ee
For the remainder $Q$ (see (\ref{2.6})) we obtain by the usual way (see (\ref{2.9})) the estimate
\be \label{2.17}
Q=\mcal{O}(e^{-\alpha T}T^{\beta-\frac 14}) .
\ee

\subsection{}

Hence, from (\ref{2.4}) by (\ref{2.6}), (\ref{2.7}), (\ref{2.13})-(\ref{2.17}) we obtain the following

\begin{lemma}
\be \label{2.18}
\begin{split}
& \Phi_1(T)=\int_T^\infty e^{-\alpha t}t^\beta Z(t){\rm d}t= \\
& =\frac{2}{\alpha}\frac{e^{-\alpha T}T^\beta}{1+\frac{1}{\alpha^2}\{\vt'(T)\}^2}\left[\cos\{\vt(T)\}-\frac{1}{\alpha}\vt'(T)\sin\{\vt(T)\}\right]+ \\
& +\frac{2}{\alpha}\sum_{2\leq n\leq\sqrt{\frac{T}{2\pi}}}\frac{e^{-\alpha T}T^\beta}{\sqrt{n}\left( 1+\frac{1}{\alpha^2}\{\vt'(T)-\ln n\}^2\right)}
\left[ \cos\{\vt(T)-T\ln n\}\frac{}{}\right.- \\
& \left. - \frac{1}{\alpha}\{\vt'(T)-\ln n\}\sin\{\vt(T)-T\ln n\}\right]+\mcal{O}(e^{-\alpha T}T^{\beta-\frac 14}) ,
\end{split}
\ee
for all sufficiently big $T>0$ and $\alpha=\frac \pi 4,\ \beta=\frac 74$.
\end{lemma}

\section{The formula for $\Psi(t_\nu)$}

\subsection{}

Since (see \cite{5}, p. 79)
\bdis
\begin{split}
& \Xi(t)=-\frac{1}{2\pi^{1/4}}\left( t^2+\frac 14\right)\left|\Gamma\left(\frac 14+i\frac t2\right)\right|Z(t), \\
& \left|\Gamma\left(\frac 14+i\frac t2\right)\right|=2^{1/4}\sqrt{2\pi}e^{-\frac{\pi t}{4}}t^{-1/4}\left\{ 1+\mcal{O}(t^{-1})\right\} ,
\end{split}
\edis
then
\be \label{3.1}
\Xi(t)=-\left(\frac{\pi}{2}\right)^{1/4}\left\{ 1+\mcal{O}(t^{-1})\right\}e^{-\frac{\pi t}{4}}t^{7/4}Z(t) .
\ee
Let
\be \label{3.2}
\Phi(T)=\int_T^\infty\Xi(t){\rm d}t=a\Phi_1(T)+\Phi_2(T)
\ee
where (see (\ref{2.18}), (\ref{3.1}))
\be \label{3.3}
\begin{split}
& \Phi_1(T)=\int_T^\infty e^{-\alpha t}t^\beta Z(t){\rm d}t,\quad a=-\left(\frac{\pi}{2}\right)^{1/4} , \\
& \Phi_2(T)=\mcal{O}\left(\int_T^\infty e^{-\alpha t}t^{\beta-1} |Z(t)|{\rm d}t\right) .
\end{split}
\ee
Since $Z(t)=\mcal{O}(t^{1/4})$ then (comp. (\ref{2.9}))
\bdis
\Phi_2(T)=\mcal{O}(e^{-\alpha T}T^{\beta-3/4}) ,
\edis
and consequently
\be \label{3.4}
\Phi(T)=a\Phi_1(T)+\mcal{O}(e^{-\alpha T}T^{\beta-3/4}) .
\ee
Hence, from (\ref{3.4}) by (\ref{2.18}), (\ref{3.3}) the formula
\be \label{3.5}
\begin{split}
& \Psi(T)=e^{\alpha T}T^{-\beta}\Phi(T)= \\
& = \frac{2a}{\alpha}\frac{1}{1+\frac{1}{\alpha^2}\{\vt'(T)\}^2}
\left[\cos\{\vt(T)\}-\frac{1}{\alpha}\vt'(T)\sin\{\vt(T)\}\right]+ \\
& +\frac{2a}{\alpha}\sum_{2\leq n\leq\sqrt{\frac{T}{2\pi}}}\frac{1}{\sqrt{n}\left( 1+\frac{1}{\alpha^2}\{\vt'(T)-\ln n\}^2\right)}
\left[\cos\{\vt(T)-T\ln n\}\frac{}{}\right.- \\
& \left. - \frac{1}{\alpha}\{\vt'(T)-\ln n\}\sin\{\vt(T)-T\ln n\} \right]+\mcal{O}(T^{-1/4})
\end{split}
\ee
follows for all sufficiently big $T>0$, and $\alpha=\frac\pi 4,\ \beta=\frac{7}{4}$.

\subsection{}

Since from (\ref{3.5}) in the case $T\to t,\ t\in [T,T+H],\ H=o(T)$, we have
\be \label{3.6}
\begin{split}
& \Psi(t)=\frac{2a}{\alpha}\frac{1}{1+\frac{1}{\alpha^2}\{\vt'(t)\}^2}\left[\cos\{\vt(t)\}-\frac{1}{\alpha}\vt'(t)\sin\{\vt(t)\}\right]+\\
& +\frac{2a}{\alpha}\sum_{2\leq n\leq\sqrt{\frac{t}{2\pi}}}\frac{1}{\sqrt{n}\left( 1+\frac{1}{\alpha^2}\{\vt'(t)-\ln n\}^2\right)}
\left[\cos\{\vt(t)-t\ln n\}-\frac{}{}\right. \\
& \left.-\frac{1}{\alpha^2}\{\vt'(t)-\ln n\}\sin\{\vt(t)-t\ln n\}\right]+\mcal{O}(t^{-1/4}) ,
\end{split}
\ee
and (see (\ref{2.2}))
\bdis
\begin{split}
& \frac{1}{1+\frac{1}{\alpha^2}\{\vt'(t)\}^2}=\frac{1}{1+\frac{1}{\alpha^2}\{\vt'(T)\}^2}+\mcal{O}
\left\{\frac{\vt'(c)\vt''(c)}{\left( 1+\frac{1}{\alpha^2}\{\vt'(t)\}^2\right)^2}H\right\}= \\
& =\frac{1}{1+\frac{1}{\alpha^2}\{\vt'(T)\}^2}+\mcal{O}\left(\frac{H}{T\ln^3T}\right) , \\
& \frac{1}{1+\frac{1}{\alpha^2}\{\vt'(t)-\ln n\}^2}=\frac{1}{1+\frac{1}{\alpha^2}\{\vt'(T)-\ln n\}^2}+\mcal{O}\left(\frac{H}{T\ln^3T}\right), \\
& \vt'(t)-\ln n=\vt'(T)-\ln n+\mcal{O}\left(\frac{H}{T}\right), \\
& \sum_{\sqrt{\frac{T}{2\pi}}\leq n\leq\sqrt{\frac{T+H}{2\pi}}}\frac{1}{\sqrt{n}}=\mcal{O}\left(\frac{H}{T^{3/4}}\right)
\end{split}
\edis
then from (\ref{3.6}) the formula
\be \label{3.7}
\begin{split}
& \Psi(t)=\frac{2a}{\alpha}\frac{1}{1+\frac{1}{\alpha^2}\{\vt'(T)\}^2}\left[\cos\{\vt(t)\}-\frac{1}{\alpha}\vt'(t)\sin\{\vt(t)\}\right]+ \\
& + \frac{2a}{\alpha}\sum_{2\leq n<P_0}\frac{1}{\sqrt{n}\left( 1+\frac{1}{\alpha^2}\{\vt'(T)-\ln n\}^2\right)}
\left[ \cos\{\vt(t)-t\ln n\}-\frac{}{}\right. \\
& \left.-\frac{1}{\alpha}\{\vt'(T)-\ln n\}\sin\{\vt(t)-t\ln n\}\right]+\mcal{O}(t^{-1/4})
\end{split}
\ee
follows for
\bdis
t\in [T,T+H],\ H=\mcal{O}\left(\frac{\sqrt{T}}{\ln T}\right),\ P_0=\sqrt{\frac{T}{2\pi}} .
\edis
Hence, from (\ref{3.7}) by (\ref{2.3}) we obtain the following

\begin{lemma}
\be \label{3.8}
\begin{split}
& \Psi(t_\nu)=\frac{2a}{\alpha}\frac{(-1)^\nu}{1+\frac{1}{\alpha^2}\{\vt'(T)\}^2}+ \\
& +\frac{2a}{\alpha}\sum_{2\leq n<P_0}\frac{1}{\sqrt{n}\left( 1+\frac{1}{\alpha^2}\{\vt'(T)-\ln n\}^2\right)}
\left\{(-1)^\nu\cos(t_\nu\ln n)+\frac{}{}\right. \\
& \left.+\frac{1}{\alpha}\{\vt'(T)-\ln n\}(-1)^\nu\sin(t_\nu\ln n)\right\}+\mcal{O}(t^{-1/4}) , \\
& t_\nu\in [T,T+H],\ H=\mcal{O}\left(\frac{\sqrt{T}}{\ln T}\right).
\end{split}
\ee
\end{lemma}

\section{Proof of the Theorem}

\subsection{}

Since (see \cite{3}, (23))
\bdis
\sum_{T\leq t_\nu\leq T+H}1=\frac{1}{\pi}H\ln P_0+\mcal{O}\left(\frac{H^2}{T}\right)
\edis
then we obtain from (\ref{3.8}) the following equalities
\be \label{4.1}
\begin{split}
& \sum_{T\leq t_\nu\leq T+H} \Psi(t_\nu)=w_1+w_2+\mcal{O}(HT^{-1/4}\ln T) , \\
& \sum_{T\leq t_\nu\leq T+H} (-1)^\nu\Psi(t_\nu)=\frac{2a}{\pi\alpha}\frac{H\ln P_0}{1+\frac{1}{\alpha^2}\{\vt'(T)\}^2}+w_3+w_4+\mcal{O}(HT^{-1/4}\ln T), \\
& H=\mcal{O}(T^{1/4+\delta})
\end{split}
\ee
for the main sums ($0<\delta$ is arbitrarily small) where
\be \label{4.2}
\begin{split}
& w_1=\sum_{2\leq n< P_0}\frac{a_n}{\sqrt{n}}\sum_{T\leq t_\nu\leq T+H}(-1)^\nu\cos(t_\nu\ln n) , \\
& w_2=\sum_{2\leq n< P_0}\frac{a_nb_n}{\sqrt{n}}\sum_{T\leq t_\nu\leq T+H}(-1)^\nu\sin(t_\nu\ln n), \\
& w_3=\sum_{2\leq n< P_0}\frac{a_n}{\sqrt{n}}\sum_{T\leq t_\nu\leq T+H}\cos(t_\nu\ln n), \\
& w_4=\sum_{2\leq n< P_0}\frac{a_nb_n}{\sqrt{n}}\sum_{T\leq t_\nu\leq T+H}\sin(t_\nu\ln n) ,
\end{split}
\ee
and (see (\ref{2.2})
\be \label{4.3}
\begin{split}
& a_n=\frac{2a}{\alpha}\frac{1}{1+\frac{1}{\alpha^2}\{\vt'(T)-\ln n\}^2} , \\
& b_n=\frac{1}{\alpha}\{\vt'(T)-\ln n\}^2=\frac{2}{\alpha}\ln\frac{P_0}{n}+\mcal{O}\left(\frac 1T\right) ,
\end{split}
\ee
($a_n$ is increasing and $a_nb_n$ is decreasing).

\subsection{}

Let us remind that we have proved (see \cite{2}, p. 38, (56), \cite{3}, (26)) for the sums

\bdis
\begin{split}
& \bar{w}_1=\sum_{2\leq n< m}\frac{1}{\sqrt{n}}\sum_{T\leq t_\nu\leq T+H}(-1)^\nu\cos(t_\nu\ln n) , \\
& \bar{w}_2=\sum_{2\leq n< m}\frac{1}{\sqrt{n}}\sum_{T\leq t_\nu\leq T+H}(-1)^\nu\sin(t_\nu\ln n), \\
& \bar{w}_3=\sum_{2\leq n< m}\frac{1}{\sqrt{n}}\sum_{T\leq t_\nu\leq T+H}\cos(t_\nu\ln n), \\
& \bar{w}_4=\sum_{2\leq n< m}\frac{1}{\sqrt{n}}\sum_{T\leq t_\nu\leq T+H}\sin(t_\nu\ln n),\quad m<P_0 ,  
\end{split}
\edis
where $\sin(t_\nu\ln n)=\cos(t_\nu\ln n-\frac \pi 2)$ the following estimates
\be \label{4.4}
\bar{w}_1,\bar{w}_2=\mcal{O}(T^{1/6+\epsilon/4}) ,
\ee
(see \cite{3}, p. 97, (2)), and
\be \label{4.5}
\bar{w}_3,\bar{w}_4=\mcal{O}(T^{1/6+\epsilon/4}) ,
\ee
(see \cite{3}, p. 48, (26) and p. 98, (5); $\Delta=\frac 16$). Since $\{ a_n\}$ and $\{ a_nb_n\}$ are monotonic sequences and (see (\ref{4.3}))
\bdis
a_n=\mcal{O}(1),\ a_nb_n=\mcal{O}(\ln P_0)=\mcal{O}(\ln T)
\edis
we obtain from (\ref{4.2}), using Abel's transformation and (\ref{4.4}), (\ref{4.5}), the following estimates
\be \label{4.6}
w_1,w_2,w_3,w_4=\mcal{O}(T^{1/6+\epsilon/2}) .
\ee

\subsection{}

Thus, from (\ref{4.1}) by (\ref{4.6}) we obtain
\be \label{4.7}
\begin{split}
& \sum_{T\leq t_\nu\leq T+H}\Psi(t_\nu)=\mcal{O}(T^{1/6+\epsilon/2}) , \\
& \sum_{T\leq t_\nu\leq T+H}(-1)^\nu\Psi(t_\nu)=\frac{2a}{\pi\alpha}\frac{H\ln P_0}{1+\frac{1}{\alpha^2}\{\vt'(T)\}^2}+\mcal{O}(T^{1/6+\epsilon/2})
\end{split}
\ee
for $H=\mcal{O}(T^{1/4+\delta})$. Since (see (\ref{2.2})
\bdis
1+\frac{1}{\alpha^2}\{\vt'(T)\}^2\sim\frac{4}{\alpha^2}\ln^2P_0 ,
\edis
we obtain from (\ref{4.7}) the following asymptotic formulae
\be \label{4.8}
\begin{split}
& \sum_{T\leq t_{2\nu}\leq T+\bar{H}}\Psi(t_{2\nu})\sim\frac{a\alpha}{\pi}\frac{\bar{H}}{\ln P_0} , \\
& \sum_{T\leq t_{2\nu+1}\leq T+\bar{H}}\Psi(t_{2\nu+1})\sim -\frac{a\alpha}{\pi}\frac{\bar{H}}{\ln P_0}, \quad \bar{H}=\frac 13 T^{1/6+\epsilon} .
\end{split}
\ee
Hence, it follows from (\ref{4.8}) that there is a zero $\omega$ of the odd order of the function $\Psi(t),\ t\in [T,T+\bar{H}]$, i.e. by (\ref{3.5}),
$\omega$ is the zero of the odd order of the function
\bdis
\Phi(t)=\int_t^\infty \Psi(\tau){\rm d}\tau .
\edis

\thanks{I would like to thank Michal Demetrian for helping me with the electronic version of this work.}

\end{document}